\documentclass[12pt]{amsart}
\usepackage{graphicx,amssymb, amsthm, amsmath, amsfonts}
\def\ds{\displaystyle}

\begin{document}

\newcommand{\Z}{\mathbb Z}
\newcommand{\Q}{\mathbb Q}
\newcommand{\N}{\mathbb N}
\newcommand{\R}{\mathbb R}
\newcommand{\C}{\mathbb C}
\newcommand{\D}{\mathbb D}
\newcommand{\A}{\mathcal A}
\newcommand{\lf}{\left \lfloor}
\newcommand{\rf}{\right \rfloor}
\newcommand{\ignore}[1]{}
\newtheorem{example}{Example}
\newtheorem{thm}{Theorem}
\newtheorem{prop}[thm]{Proposition}
\newtheorem{lemma}[thm]{Lemma}
\newtheorem{cor}[thm]{Corollary}
\newtheorem{conj}[thm]{Conjecture}
\newtheorem{definition}[thm]{Definition}
\newtheorem{question}[thm]{Question}

\allowdisplaybreaks

\title{Cycles in Austrian Solitaire}             
    \author{Philip P. Mummert}                              
 \address{Department of Mathematics\\
Purdue University\\
     150 N University St\\
West Lafayette, IN 47907}                                          
  \email{pmummert@purdue.edu}                    %
    \date{\today}                                       
    \keywords{integer partitions, discrete dynamical systems, state space, asset inventory, recreational mathematics, Farey sequence, Raney's Lemma, maximally even configuration, Bulgarian Solitaire, periodic cycles}           %\keywords                      
                       
    \begin{abstract} Austrian Solitaire is a variation of Bulgarian Solitaire. It may be described as a card game, a method of asset inventory management, or a discrete dynamical system on integer partitions. We prove that the limit cycles in Austrian Solitaire do not depend on the initial configuration; in other words, each state space is connected. We show that a full Farey sequence completely characterizes these unique (and balanced) cycles.
\end{abstract} 
\maketitle

\section{Introduction} The recreational mathematics puzzle known as ``Bulgarian Solitaire'' was presented to the readership of {\it Scientific American} by Martin Gardner in 1983 \cite{G}.  Begin with a deck of $n$ cards distributed into piles. To play Bulgarian Solitaire, draw one card from each pile and form one new pile with all the drawn cards.  Repeat.  Easily explained as a procedure of pile-making from a deck of cards, one might also describe it in more sophisticated terms as an iterated operation on integer partitions. For a detailed account of Bulgarian Solitaire's mathematical history, see \cite{H}.    See \cite{Br} for solutions to counting and characterizing cycles of any deck size $n$.  If $n$ is a triangular number, there is a unique fixed state which is the limiting outcome realized by all possible initial pile configurations.  For example, with $n=6$ cards, the game always leads to a pile of 3, a pile of 2, and a pile of 1.  But not every deck size has a unique cycle; $4+2+2$ and $3+3+2$ lie in different cycles for $n=8$.

In their 1985 exposition on the Bulgarian Solitaire solution, Akin and Davis \cite{AD} introduce  ``Austrian Solitaire,'' a variant requiring an additional parameter.  Begin with a deck of $n$ cards distributed into piles. To play Austrian Solitaire, draw one card from each pile and form as many new piles of $L$ cards as possible, leaving the remainder in hand for the next turn.  Repeat. What is the end result of playing such a game?  Might the outcome depend on the initial state?   In this paper we prove the conjecture first presented by Akin and Davis and reiterated by da Silva, Hopkins, and Sellers \cite{DHS}: there is only one cycle possible for each fixed $n$ and $L$.  Other related combinatorial questions have been considered previously by others \cite{Ba}, \cite{DHS}.

To make some sense of the ``Austrian'' moniker,   we reimagine the process using the language of capital and reinvestment in the spirit of Austrian economic theory,  and thus lend our affirmation to the notion expressed in \cite{AD} when they declare this ``economic interpretation shows that anything can become applied mathematics.''
  Suppose a company maintains an asset inventory of items (e.g., a collection of laptops, a fleet of cars, etc.) that are used in operations.  These items lose value and each one needs replaced every $L$ years.  The company would like to implement an easy replacement cycle that keeps the total value of the asset inventory constant.  So each year the budget sets aside an amount equal to the prior year's total depreciation, and as many new items are purchased as possible.  Any leftover funds roll over for future use.  Now over time, we hope that such a procedure will eventually lead to a rather balanced asset inventory with roughly the same number of items in use each year, and the age of those items uniformly spread out, but can one prove that this sort of scenario is realized in all cases?

Translating from piles of cards or asset management into mathematical notation, we fix natural numbers $L$ and $n$, and say $\lambda=(\lambda_0 ; \lambda_1, \ldots, \lambda_t)$ is an {\it Austrian partition} of the integer $n$ with {\it bank capacity} $L$ if $\lambda_0, \ldots, \lambda_t$ are integers such that $\lambda_0+\ldots+\lambda_t=n$, $0\leq \lambda_0 <L$, and $L\geq \lambda_1\geq \lambda_2 \geq \ldots \lambda_t \geq 1$.  We call $\lambda_0$ the {\it bank deposit}.  
Let $f_m(\lambda)$ denote the frequency of the parts of size $m$ (not including the bank deposit) in the partition $\lambda$, i.e. \[f_m(\lambda)=\#\{\lambda_i=m: i\neq 0\},\]
If $t(\lambda)$ denotes the number of (non-bank) parts of $\lambda$,
 {\it Austrian Solitaire} is played by iterating the function  \[ A(\lambda)= (\lambda_0+t(\lambda)-kL; \ \overbrace{L, \ldots, L}^{k \textnormal{ times}}, \lambda_1-1, \lambda_2-1, \ldots, \lambda_{s}-1);\]
 $k$ is the greatest integer less than or equal to $\ds \frac{\lambda_0+t(\lambda)}{L}$ and $\ds s=t(\lambda)-f_1(\lambda)$.  If our initial partition of $n$ is not Austrian (i.e., $\lambda_0\ge L$ or $\lambda_i>L$ for some $i>0$), after a few iterations we obtain one that is Austrian, and so we may restrict our attention to these partitions.  Let $\A (n,L)$ denote the set of Austrian partitions of $n$ with bank capacity $L$.  For any $\lambda\in \mathcal A(n,L)$,  \[n=\lambda_0+\sum_{m=1}^L m\cdot f_m(\lambda).\] We will use $\lfloor \cdot \rfloor$ and $\lceil \cdot \rceil$ to denote the integer floor and ceiling functions, respectively.   
For each $\lambda \in \mathcal A(n,L)$, the partition $A(\lambda)\in \mathcal A(n,L)$ has  
\[\ds f_m(A(\lambda)) = f_{m+1}(\lambda)\ \ \text{ for } 1\leq m\leq L-1, \] \[\ds f_L(A(\lambda)) = \left \lfloor \frac{1}{L}\left (\lambda_0+ \sum_{m=1}^L f_m(\lambda) \right )\right \rfloor,\]
\[\text{and bank deposit } \ [A(\lambda)]_0=\ds \left (\lambda_0+\sum_{m=1}^L f_m(\lambda) \right )\textrm{mod } L.\]

Because the set of states is finite for each fixed pair of $n$ and $L$ values, every initial state must iterate to a cycle.  We will show the cycle is unique and consequently, we have a totally connected state diagram. In other words, the cycle reached under iteration of $A$ depends only on $n$ and $L$ and does not depend on the choice of initial partition.  
In terms of asset inventory replacement, that means  the long-term cycle is uniquely determined by only two parameters: the (constant) total value and the item lifespan; it does not depend on the initial configuration of the inventory.
 Furthermore, the limit cycle is characterized by means of what we will call a ``full'' Farey sequence: let $a_0, a_1, a_2, \ldots$ be the sequence of non-negative fractions with denominator no greater than $L$ listed in non-decreasing order.  Our sequence is a modification of the standard Farey sequence, because we will list all equivalent fractions as separate entries; the full Farey sequence includes all possible numerator and denominator pairs, and not just those in lowest terms. 
For example, the full Farey sequence for $L=5$ is
\[\frac{0}{1},\frac{0}{2},\frac{0}{3},\frac{0}{4},\frac{0}{5}, \frac{1}{5}, \frac{1}{4}, \frac{1}{3}, \frac{2}{5}, \frac{1}{2}, \frac{2}{4}, \frac{3}{5}, \frac{2}{3}, \frac{3}{4}, \frac{4}{5}, \frac{1}{1}, \frac{2}{2},\frac{3}{3}, \frac{4}{4}, \frac{5}{5}, \frac{6}{5}, \frac{5}{4},  \frac{4}{3}, \frac{7}{5}, \cdots \]

 Austrian partitions belonging to a cycle must distribute part sizes in a uniquely balanced (or maximally-even) way, and all cases of this distribution are encoded by the full Farey sequence:

\begin{thm} \label{theorem} Fix $n\ge 0$ and $L\ge 1$.  Let $a_0, a_1, a_2, \ldots$ be the ``full'' Farey sequence with denominators no greater than $L$.   Write $a_n$ as $\dfrac{q}{p}$, in simplest terms.   With a deck of $n$ cards and $L$ as the size of new piles, any initial state iterates to the same cycle.  This unique cycle has period $p$.  The partition with smallest bank deposit within this cycle satisfies $\ds f_m(\lambda)=\left \lceil \frac{q}{p}(L-m+1)  \right \rceil -\left \lceil \frac{q}{p}(L-m) \right \rceil$ for $1\le m \le L$, and has bank deposit equal to the number of times $a_k=\dfrac{q}{p}$ for $k<n$. 
\end{thm}

\begin{example} The full Farey sequence for $L=5$ has $a_{22}=\dfrac{4}{3}$.  So the Austrian cycle for $n=22$ and $L=5$ has period $3$ and the Austrian partition with smallest bank within this 3-cycle  is $(0; 5,5,4,3,2,2,1)$.
Every $\lambda \in \mathcal A(22,5)$ will reach this unique cycle under iteration of $A$. \end{example}

\begin{example}  The full Farey sequence for $L=14$ has $a_{139}=\dfrac{5}{4}$.  So the Austrian cycle for $n=139$ and $L=14$ has period $4$ and the Austrian partition with smallest bank within this 4-cycle  is
\[(2; 14,14,13,12,11,10,10,9,8,7,6,6,5,4,3,2,2,1).\]  The bank deposit is $2$ because  $a_{137}=a_{138}=a_{139}=\dfrac{5}{4}=\dfrac{10}{8}=\dfrac{15}{12}$.
Every $\lambda \in \mathcal A(139,14)$ will reach this unique cycle under iteration of $A$.\end{example}

Why is the theorem true?  As a result of the simplicity of the Austrian Solitaire action, all candidates for cyclic behavior of period $p$ can be concisely described  with a $p$-tuple of non-negative integers: the cyclic progression of frequencies of any part size (which is the same for every part size).  Very few of these cyclic possibilities can be realized;  if we regard the bank deposit as merely the difference between $n$ and $\ds \sum_{m=1}^L m\cdot f_m(\lambda)$ for each $\lambda$ in the cycle, the requirement that this deposit must always be non-negative and strictly less than $L$ turns out to be very restrictive.  In fact, there is only one such cycle for any fixed $n$ and $L$.

Here is the proof outline:  First, we give a simple way to catalog any potential periodic cycles, which conjugates our dynamical system to a simple shift map.  In considering the distribution of part frequencies within a periodic cycle, we focus on part distributions that are, in some sense, evenly spread out or ``balanced.''   By using a maximal state within the balanced cycle as our plumb line, we show that imbalanced solutions cannot be realized due to the bank constraint.  This particular step is the key component of the proof and relies on an argument conceptually similar to the proof of Raney's Lemma found in \cite{GKP}.
As the final step toward obtaining uniqueness, we show two distinct balanced cycles can't be realized by the same set of $n$ and $L$ values.  We conclude by summarizing the result as it relates to our full Farey sequence.

Before proceding to the proof of Theorem \ref{theorem}, we provide a corollary that gives a simple formula for the index of a fraction in the full Farey sequence.  Since $\ds n=\lambda_0+\sum_{m=1}^Lm\cdot f_m(\lambda)$, the partition described in Theorem \ref{theorem} immediately gives the following consequence:
\begin{cor} \label{corollary} The fraction $\dfrac{q}{p}$ (with $1\le p\le L$)  first appears as $a_N$ in the full Farey sequence (with denominators no greater than $L$) for $N=\displaystyle \sum_{k=1}^L\left \lceil \frac{kq}{p}\right \rceil$.  In fact, $a_n=\dfrac{q}{p}$ where $N\le n < N+\left \lfloor\dfrac{L}{p}\right \rfloor $.
\end{cor}

\section{Proof of the Theorem}

%Consider $p$ non-negative integers, $\beta_0, \beta_1, \beta_2, \ldots, \beta_{p-1}$.  We extend to a bi-infinite sequence $\beta$ by setting $\beta_i=\beta_{i \textrm{ mod }p}$ for all integers $i$.  

 We use exponent notation on functions to denote iteration.  Assume $\lambda \in \mathcal A(n, L)$ is part of a periodic cycle, with $A^p(\lambda)=\lambda$.  
The meaning of $A^{-1}$ is well defined within this periodic cycle.  Define $\beta$, a bi-infinite sequence of non-negative integers by $\beta_{i}=f_L(A^{-i}(\lambda) )$ for each integer $i$.  We write $\phi$ for this mapping, $\beta=\phi(\lambda)$.
Let $\sigma$ be the shift map given by $(\sigma \beta)_i=\beta_{i-1}$.  The effect of $A$ on pile sizes within a cycle gives us \[\sigma (\phi (\lambda)) = \phi (A(\lambda))\]
and $\beta=\phi(\lambda)$ must also be periodic: $\beta_i=\beta_{i \text{ mod } p}$.
For fixed values of $n$ and $L$, $\phi$ is injective on the set of periodic partitions in $\mathcal A(n,L)$.  In other words, $\beta=\phi(\lambda)$ completely determines $\lambda$, since \[f_m(\lambda) = f_L(A^{m-L}\lambda)=\beta_{L-m}\] for $1\le m \le L$, and the bank deposit, $\ds \lambda_0=n-\sum_{m=1}^L m\cdot f_m(\lambda)$.

%%%%%
\begin{example} If the sequence $\beta=\ldots 10010101102000\ldots$ (with $p=14$) were equal to 
$\phi(\lambda)$ for some $\lambda \in \mathcal A(n,20)$ belonging to a cycle, then \[\lambda=(\lambda_0;20,17,15,13,12,10,10,6,3,1)\] with $\lambda_0=n-107$.  In fact, we will show this is not possible.
\end{example}

%For notational convenience, partial sums will be denoted as \[\beta_m^k=\ds \sum_{i=0}^{k-1}\beta_{m+i}\] where $k$ and $m$ are integers with $k\ge 1$.
The sum of partial sums gives the total of the non-bank parts:
\begin{eqnarray*}
 \sum_{k=1}^L   \sum_{i=0}^{k-1}\beta_i&=&\sum_{i=0}^{L-1}   \sum_{k=i+1}^{L}\beta_i = \sum_{m=1}^{L}m\cdot \beta_{L-m}\\
&=& \sum_{m=1}^{L}m\cdot f_L(A^{m-L}\lambda)=\sum_{m=1}^{L}m\cdot f_m(\lambda)=n-\lambda_0.
\end{eqnarray*}

We call a bi-infinite sequence $\beta$ ``balanced'' if for some $j$, $q$, and $p$, $\gamma=\sigma^j\beta$ satisfies 
\begin{eqnarray}
\gamma_i=\left \lceil \frac{q}{p}(i+1)  \right \rceil -\left \lceil \frac{q}{p}i \right \rceil
\end{eqnarray} for all integers $i$.  
We call this particular iterate, $\gamma=\sigma^j \beta$, ``maximal.''   Note that $\beta$ has period $p$ and we can assume $q$ and $p$ are relatively prime.   Our notion of ``balanced'' is equivalent to what has been called ``maximally even'' in some other contexts \cite{DK}.%, \cite{DK2}.

For each integer $m$, the telescoping sum $\ds \sum_{i=0}^{k-1} \gamma_{m+i}$ is equal to $ \ds \left \lceil \frac{q}{p}(m+k)\right  \rceil - \left \lceil \frac{q}{p}m\right  \rceil$, and, in particular, $\ds \sum_{i=0}^{k-1} \gamma_{i}= \ds \left \lceil \frac{q}{p}k\right  \rceil$.
In general,  $\lceil a +b \rceil \leq \lceil a \rceil + \lceil b \rceil$, so  $\ds \left \lceil \frac{q}{p}(m+k)\right  \rceil - \left \lceil \frac{q}{p}m\right  \rceil \le \left \lceil \frac{q}{p}k \right \rceil$, thus, $\ds \sum_{i=0}^{k-1} \gamma_{m+i} \le  \sum_{i=0}^{k-1} \gamma_{i}$ for any $m$ and $k$, which means $\gamma$ is maximal in partial sums among its iterated shifts.

%!!!!!!!!!!!!!!!!!!!!!!!!!!!!!!!!!!!!!!!!!!!!!!!!!!!!!!!!!!!!!!!!!!!!!!!!!!!!!!!!!!!!!!!!!!!!!!!!!!!

%Note $\gamma$ has least period $p/\textrm{gcd}(p,q)$.  So we can now assume that our $\beta$ is balanced, that $p$ is the least period, and importantly, that $q$ is relatively prime to $p$.

\begin{prop}\label{gamma}%stated more precisely
Assuming $q$ and $p$ are relatively prime,  the maximal balanced sequence, $\gamma$, given by \[ \gamma_i=\left \lceil \frac{q}{p}(i+1)  \right \rceil -\left \lceil \frac{q}{p}i \right \rceil\] is realized as $\phi(\lambda)$ for some periodic $\lambda \in \mathcal A(n,L)$ if and only if the bank deposit $\lambda_0$ is less than $\ds \left \lfloor \frac{L}{p}\right \rfloor$.  
\end{prop}

%Notice that since $\lambda_0\ge 0$, we must have $L\ge p$.

\begin{proof}
Consider the plot of the partial sums, $\ds \sum_{i=0}^{k-1} \gamma_i$ as a function of $k$.  Look at the least upper bounding line of slope $\ds \frac{q}{p}$, and identify one of the intersection points as the start of an iterated shift of $\gamma$.  In other words, if $\ds B=\min \{b\in \R: \sum_{i=0}^{k-1} \gamma_i\le b+\frac{q}{p}k, \textrm{ for each } k\}$, then $\ds \sum_{i=0}^{j-1} \gamma_i= B+\frac{q}{p}j$ for some particular $j\ge 1$.
If $\ds \sum_{i=0}^{k-1}\gamma_{j+i}>\frac{q}{p}k$, then $\ds \sum_{i=0}^{j+k-1}\gamma_i=\sum_{i=0}^{j-1}\gamma_i+\sum_{i=0}^{k-1}\gamma_{j+i}>B+\frac{q}{p}\left (j +k \right )$, a contradiction.
Hence $\ds \sum_{i=0}^{k-1} \gamma_{j+i}\le \frac{q}{p}k$, and since it is integer-valued,   $\ds \sum_{i=0}^{k-1}\gamma_{j+i}\le \left \lfloor \frac{q}{p}k \right \rfloor$.
In general, $\lceil a +b \rceil \geq \lceil a \rceil + \lfloor b \rfloor$, so $\ds \sum_{i=0}^{k-1}\gamma_{m+i}= \left \lceil \frac{q}{p}(m+k)\right  \rceil - \left \lceil \frac{q}{p}m\right  \rceil \ge \left \lfloor \frac{q}{p}k \right \rfloor$ for any $m$ and $k$.  Hence, $\ds \sum_{i=0}^{k-1}\gamma_{j+i}=\left \lfloor \frac{q}{p}k \right \rfloor$ and, in summary, for any integers $m$ and $k$ with $k\ge 1$: \[ \left \lfloor \frac{q}{p}k \right \rfloor=\sum_{i=0}^{k-1} \gamma_{j+i} \le \sum_{i=0}^{k-1} \gamma_{m+i} \le \sum_{i=0}^{k-1} \gamma_{i} =\left \lceil \frac{q}{p}k \right \rceil\]
Furthermore,
\[ \sum_{k=1}^L \sum_{i=0}^{k-1} \left (\gamma_{i}- \gamma_{j+i}\right )=\sum_{k=1}^L \left (\left \lceil \frac{q}{p}k \right \rceil-\left \lfloor \frac{q}{p}k \right \rfloor \right )=L-\left \lfloor \frac{L}{p}\right \rfloor \]
where the last equality is given by the fact that $q$ is relatively prime to $p$; the floor and ceiling are equal if and only if $k$ is a multiple of $p$. 

Now if $\phi(\lambda)=\gamma$ for some $\lambda \in \mathcal A(n,L)$ belonging to a periodic cycle, then the bank deposits within the cycle must be related as follows:

\begin{eqnarray*}\lambda_0= \left (n-\sum_{k=1}^L \sum_{i=0}^{k-1}\gamma_{i}\right )&\le&\left ( n-\sum_{k=1}^L \sum_{i=0}^{k-1} \gamma_{m+i}\right )\\ &\le& \left (n-\sum_{k=1}^L \sum_{i=0}^{k-1}\gamma_{j+i}\right ) =\lambda_0+L-\left \lfloor \frac{L}{p}\right \rfloor
\end{eqnarray*}
for any integer $m$.
The constraint on bank deposits requires  $\ds \lambda_0+L-\left \lfloor \frac{L}{p}\right \rfloor<L$, so $\gamma$ can be realized by a partition $\lambda$ if and only if $\ds 0\le \lambda_0< \left \lfloor \frac{L}{p}\right \rfloor$.
\end{proof}

%***************************************************************

Next, we show that an Austrian partition in a periodic cycle cannot correspond to an imbalanced sequence.  We compare an appropriate shift iterate to the maximal balanced solution $\gamma$, and find that any differences (i.e. imbalances) imply the largest and smallest bank deposits within the cycle must differ by at least $L$, which contradicts the constraint on the bank capacity.

\begin{figure}[h]
\centering
\includegraphics[scale=0.6]{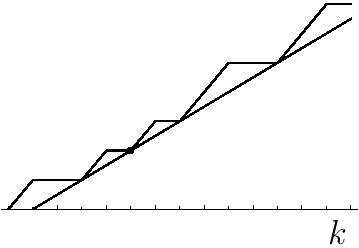}

\caption{The plot of $\ds \sum_{i=0}^{k-1}\beta_i$  and its greatest lower bounding line of slope $\dfrac{7}{14}$ for $\beta=\ldots 10010101100110\ldots$  .}

\end{figure}

\begin{prop}\label{balance} Suppose $\lambda \in \mathcal A(n, L)$ belongs to a periodic cycle.  If $\beta=\phi(\lambda)$, then $\beta$ is balanced.
\end{prop}

\begin{proof}

Suppose $\beta$ has period $p$, and let $\ds q=\sum_{i=0}^{p-1}\beta_i$.   
Consider the plot of the partial sums, $\ds \sum_{i=0}^{k-1} \beta_i$ as a function of $k$.  Consider the greatest lower bounding line of slope $\ds \frac{q}{p}$, and identify one of the intersection points as the start of $\omega$, an iterated shift of $\beta$.  In other words, if $\ds B=\max \{b\in \R: \sum_{i=0}^{k-1}\beta_i\ge b+\frac{q}{p}k, \textrm{ for each } k\}$, then $\ds \sum_{i=0}^{j-1}\beta_i= B+\frac{q}{p}j$ for some $j$, and we let $\omega=\sigma^{-j}\beta$. 
If $\ds \sum_{i=0}^{k-1} \omega_i=\sum_{i=0}^{k-1}\beta_{j+i}<\frac{q}{p}k$, then $\ds\sum_{i=0}^{j+k-1} \beta_i=\sum_{i=0}^{j-1}\beta_i+\sum_{i=0}^{k-1}\beta_{j+i}<B+\frac{q}{p}\left (j +k \right )$, a contradiction.

Hence $\ds \sum_{i=0}^{k-1}\omega_i\ge \frac{q}{p}k$, and since it is integer-valued,  $\ds \sum_{i=0}^{k-1}\omega_i\ge \left \lceil \frac{q}{p}k \right \rceil$, so $\ds\sum_{i=0}^{k-1}\omega_i\ge \sum_{i=0}^{k-1}\gamma_i$ for each $k\ge 1$, where $\gamma$ is the maximal balanced solution for $p$ and $q$.

Now suppose $\omega \neq \gamma$.  Let $\ds J=\max\{1\le j \le p: \sum_{i=0}^{j-1}\omega_i> \sum_{i=0}^{j-1}\gamma_i\}$.   It is important to see that $J<p$ since  $\ds \sum_{i=0}^{p-1} \omega_i=\sum_{i=0}^{p-1}\gamma_i=q$. 
%In general, $\ds \omega_j=\sum_{i=0}^{j}\omega_i-\sum_{i=0}^{j-1}\omega_i$ and $\ds \gamma_j=\sum_{i=0}^{j}\gamma_i-\sum_{i=0}^{j-1}\gamma_i$ for any $j\ge 1$. 
For $J+1\le j\le p$, we have $\ds \sum_{i=0}^{j-1}\omega_i= \sum_{i=0}^{j-1}\gamma_i$, and so it follows that $\omega_j=\gamma_j$ for $J+1\le j \le p-1$.
And since $\ds \sum_{i=0}^{J}\omega_i= \sum_{i=0}^{J}\gamma_i$, but $\ds \sum_{i=0}^{J-1}\omega_i> \sum_{i=0}^{J-1}\gamma_i$, we see that  $ \gamma_{J}\ge \omega_{J}+1$.
Recall that $\ds \sum_{i=0}^{k-1}\gamma_i \ge \sum_{i=0}^{k-1}\gamma_{J+i}$ for any $k\ge 1$ by the maximal property of $\gamma$.
In what follows we adopt the convention that in the cases where our lower limit of summation exceeeds the upper limit, the empty sum is equal to 0.

\begin{eqnarray*}
&&\sum_{k=1}^L  \sum_{i=0}^{k-1}\left (\gamma_{J+i} - \omega_{J+i} \right )\\
&=& \sum_{k=1}^L \left ( \gamma_{J+0} - \omega_{J+0} \right ) + \sum_{k=2}^L  \sum_{i=1}^{k-1}\left (\gamma_{J+i} - \omega_{J+i} \right )\\
&=& L \cdot \left ( \gamma_J - \omega_J \right ) + \sum_{k=2}^L  \sum_{i=0}^{k-2}\left (\gamma_{J+1+i} - \omega_{J+1+i} \right )\\
 &\ge& L\cdot 1 + \sum_{k=1}^{L-1}\sum_{i=0}^{k-1}  \left (\gamma_{J+1+i}- \omega_{J+1+i}\right )\\
&=& L + \sum_{k=1}^{p-J-1} \sum_{i=0}^{k-1} \left (  \gamma_{J+1+i} - \omega_{J+1+i}\right ) + \sum_{k=p-J}^{L-1}  \sum_{i=0}^{k-1}\left (\gamma_{J+1+i} - \omega_{J+1+i} \right )\\
%&=& L + 0 +  \sum_{k=M+1}^{L-1} \left ( \sum_{i=0}^{M-1}\gamma_{J+1+i} +\sum_{i=0}^{k-M-1}\gamma_i - \sum_{i=0}^{M-1}\omega_{J+1+i} - \sum_{i=0}^{k-M-1}\omega_{i}\right )\\
&=& L + 0 + \sum_{k=p-J}^{L-1} \sum_{i=0}^{p-J-2}\left (\gamma_{J+1+i} - \omega_{J+1+i} \right )+\sum_{k=p-J}^{L-1}  \sum_{i=p-J-1}^{k-1}\left (\gamma_{J+1+i} - \omega_{J+1+i} \right ) \\
&=& L + \sum_{k=p-J}^{L-1} \sum_{i=J+1}^{p-1}\left (\gamma_{i} - \omega_{i} \right )+\sum_{k=1}^{L-p+J}  \sum_{i=p}^{p+k-1}\left (\gamma_{i} - \omega_{i} \right ) \\
&=& L + 0 +\sum_{k=1}^{L-p+J}  \sum_{i=p}^{p+k-1}\left (\gamma_{i} - \omega_{i} \right ) = L +\sum_{k=1}^{L-p+J} \sum_{i=0}^{k-1} \left (\gamma_i - \omega_i\right )
\end{eqnarray*}
Hence,
\begin{eqnarray*}
&& \sum_{k=1}^L \sum_{i=0}^{k-1} \left ( \omega_i - \omega_{J+i} \right )\\
&=&\sum_{k=1}^L \sum_{i=0}^{k-1}\left ( \omega_i - \gamma_i \right )+\sum_{k=1}^L  \sum_{i=0}^{k-1}\left ( \gamma_i - \gamma_{J+i} \right )+\sum_{k=1}^L \sum_{i=0}^{k-1} \left (\gamma_{J+i}- \omega_{J+i} \right )\\
 &\ge& \sum_{k=1}^L \sum_{i=0}^{k-1} \left (\omega_i - \gamma_i \right )+ 0+ \left [L+ \sum_{k=1}^{L-p+J}  \sum_{i=0}^{k-1} \left ( \gamma_i - \omega_i\right ) \right ]\\
 &\ge& L+\sum_{k=M}^{L}\sum_{i=0}^{k-1}  \left (\omega_i - \gamma_i\right ) \ \ \ \ \ \   \text{where $M=L-p+J+1\le L$}\\
 &\ge& L\\
\end{eqnarray*}
Therefore, assuming that $\omega \neq \gamma$ results in $\ds \sum_{k=1}^L \sum_{i=0}^{k-1} \left ( \omega_i - \omega_{J+i} \right )\ge L$. 
Relating this inequality back to $n$ and bank deposits gives 
\begin{eqnarray*}
L\le  \sum_{k=1}^L \sum_{i=0}^{k-1} \omega_i - \sum_{k=1}^L \sum_{i=0}^{k-1} \omega_{J+i}&=&(n-[A^{-j}\beta]_0)-(n-[A^{J-j}\beta]_0)\\&=&[A^{J-j}\beta]_0-[A^{-j}\beta]_0 < L,\end{eqnarray*}
a contradiction.  Thus, $\omega=\gamma$ and  $\beta$ is balanced.

\end{proof}

%%%%%%%%%%%%%%%%%%%%%%%%%%%%%%%%%%%%%%%%%%%%%%%
Now that we have eliminated imbalanced solutions as a possibility, might it be that case that more than one balanced solution could be realized for the same $n$?  The answer is no:

\begin{prop}\label{unique} Suppose $\lambda \in \mathcal A(n,L)$ has period $p$ and $\lambda' \in \mathcal A(n',L)$ has period $p'$.  If $\ds \frac{q'}{p'}>\frac{q}{p}$, where $\ds q=\sum_{i=0}^{p-1}[\phi(\lambda)]_i$ and $\ds q'=\sum_{i=0}^{p'-1}[\phi(\lambda')]_i$, then $\ds n'>n$.
\end{prop}

\begin{proof}
By replacing $\lambda$ with some iterate, we may assume that $\phi(\lambda)=\gamma$ is maximal and balanced, and likewise for $\lambda'$.
Suppose $\ds \frac{q'}{p'}>\frac{q}{p}$. 

\begin{eqnarray*}
 \sum_{k=1}^L \sum_{i=0}^{k-1}\left ( \gamma'_i -\gamma_i \right )&=& \sum_{k=1}^L \left (\left \lceil \frac{q'}{p'}k \right \rceil-\left \lceil \frac{q}{p}k \right \rceil \right )
\end{eqnarray*}
is a sum of non-negative terms.  Since $\ds \frac{p}{p'}q'>q$, we see that when $k$ is a multiple of $p$, i.e. $k=\ell p$,  the summand is strictly positive:
 $\ds \left \lceil \frac{q'}{p'}k \right \rceil-\left \lceil \frac{q}{p}k \right \rceil =\left \lceil \frac{p}{p'}q'\ell \right \rceil - q\ell\ge 1$.
 So $\ds \sum_{k=1}^L \sum_{i=0}^{k-1}\left ( \gamma'_i -\gamma_i \right )\ge \left \lfloor \frac{L}{p} \right \rfloor$. 
By Proposition \ref{gamma}, and because the bank deposit for $\lambda'$ is non-negative, 
\begin{eqnarray*}
 n' &\ge&\sum_{k=1}^L \sum_{i=0}^{k-1}\gamma'_i \  \ge \ \sum_{k=1}^L \sum_{i=0}^{k-1} \gamma_i + \left \lfloor \frac{L}{p} \right \rfloor \\
  &=& n-\lambda_0 + \left \lfloor \frac{L}{p} \right \rfloor > n.
  \end{eqnarray*}
\end{proof}

\noindent We conclude by connecting our results back to the full Farey sequence.

\begin{proof}[Proof of Theorem 1]
 Fix $L\ge 1$ and $n\ge 0$.  There are only finitely many states possible, so a periodic cycle must exist.  By Proposition \ref{balance}, this cycle is balanced, meaning it has some iterate $\lambda \in \mathcal A(n,L)$ with $\phi(\lambda)=\gamma$, for a maximal balanced $\gamma$ as defined in Equation (1) with some fraction  $\dfrac{q}{p}$, where $q$ and $p$ are relatively prime.    Now if there is another periodic cycle for these values of $n$ and $L$, we must have $\gamma=\phi(\mu)$ for some $\mu$ in the cycle, by the contrapositive of Proposition \ref{unique} and the fact that $\dfrac{q}{p}$ defines $\gamma$.  Hence, all the non-bank parts of $\lambda$ and $\mu$ are identical, and since we assumed $n$ is the same for both, the bank deposits agree as well.  In other words, $\lambda$ must equal $\mu$ and the periodic cycle is unique.

Let us catalog our periodic Austrian cycles.   Fix $L\ge 1$ and consider the sequence (as $n\ge 0$ increases) of fractions $\dfrac{q}{p}$ which characterize the periodic cycles.  Which fractions occur, and in what order?  For $n\in \left \{0,1,\ldots, L-1\right \} $, it is apparent that the only periodic cycle is the fixed point consisting simply of bank $\lambda_0=n$, since no new piles can be formed.  So we begin with the fraction $\dfrac{0}{1}$ for each $0\le n \le L-1$.  Since $\lambda_0\ge 0$, Proposition \ref{gamma} indicates that $\dfrac{q}{p}$ is realized for some $\lambda$ if and only if $p\le L$.  Hence, our list must contain all the fractions (in lowest terms) with $p\le L$.  And each of these fractions must be repeated precisely $\left \lfloor\dfrac{L}{p}\right \rfloor $ times, once for each possible value of $\lambda_0$,  according to Proposition \ref{gamma}.  One way to enumerate these possibilities is by listing all the equivalent fractions $\ds \frac{q}{p}=\frac{2q}{2p}=\frac{3q}{3p}=\ldots=\frac{kq}{kp}$, where $kp\le L$, i.e., $\ds k\le  \left \lfloor \frac{L}{p} \right \rfloor$.  Furthermore, the increasing relation between $\dfrac{q}{p}$ and $n$ in Proposition \ref{unique} means our fractions must be listed in non-decreasing order.  Hence, our full Farey sequence, which lists all numerator and denominator pairs $\dfrac{q}{p}$, with $q\ge 0$ and $1\le p \le L$, assigns a unique cycle for each  $n$ exactly in the manner we have described here.

%Fix $n$ and $L$.    If $\lambda \in \mathcal A(n,L)$ belongs to a periodic cycle, by Proposition  \ref{balance}, $\beta=\phi(\lambda)$ is a balanced sequence generated by some $\dfrac{q}{p}$. Recalling that $\phi$ is  injective on the set of periodic partitions in $\mathcal A(n,L)$, by Proposition \ref{unique}, there does not exist another $\lambda' \in \mathcal A(n,L)$ that is periodic with $\phi(\lambda')$ generated by some other $\dfrac{q'}{p'}\neq \dfrac{q}{p}$.  Because $n=\lambda_0+\ds \sum_{k=1}^L \beta_0^k$, if $\lambda \in \mathcal A(n,L)$ has least period $p$, by Proposition \ref{gamma}, there are exactly $\ds  \left \lfloor \frac{L}{p} \right \rfloor$ possible values for $\lambda_0$, and hence, for $n$.Hence, our full Farey sequence, which lists all numerator and denominator pairs $\dfrac{q}{p}$, with $q\ge 0$ and $1\le p \le L$, assigns a unique cycle for each  $n$. 

\end{proof}

\section{Acknowledgements}  The author's interest in this problem began with a research project investigation by Taylor University undergraduates Daniel Kasper and Rachel DeMeo in 2010.  
Having now proven the conjecture of Akin and Davis, we feel some obligation to fulfill their concluding promise in \cite{AD} by asserting that our readers will now be able to say, when asked about business cycles, ``Why, it’s all in the cards!''


\begin{thebibliography}{99}

\bibitem[AD]{AD}
E Akin and M Davis, {\it Bulgarian Solitaire}, Amer. Math. Monthly {\bf 92} (1985), no. 4, 237--250.

\bibitem[Ba]{Ba}
K Bastola, {\it Enumeration of Austrian Solitaire}, Baccalaureate Honors Thesis, Saint Peter’s University Digital Repository, Jersey City NJ, 2012.

\bibitem[Br]{Br}
J Brandt, {\it Cycles of Partitions}, Proc. Amer. Math. Soc {\bf 85} (1982), no. 3, 483--486.


%\bibitem[DK1]{DK1} 
%Douthett, J., Krantz, R. {\it Maximally even sets and configurations: common threads in mathematics, physics, and music}, J. Comb. Optim. 14, 385–410 (2007).

\bibitem[DHS]{DHS}
R da Silva, B Hopkins, and J Sellers, {\it Garden of Eden States in Austrian Solitaire}, European J. Combin. {\bf 83} (2020), [103023].


\bibitem[DK]{DK} 
J Douthett and R Krantz, {\it Dinner Tables and Concentric Circles: A Harmony of Mathematics, Music, and Physics}, Col. Math. J. {\bf 39} (2008), no. 3, 203--211.

\bibitem[G]{G}
M Gardner, {\it Mathematical Games: Tasks you cannot help finishing no matter how hard you try
to block finishing them}, Sci. Amer. {\bf 249} (1983), no. 2, 12--21.


\bibitem[GKP]{GKP}
R Graham, D Knuth, and O Patashnik, {\it Concrete Mathematics: A Foundation for Computer Science}, Addison Wesley, 1994, p. 360. 

\bibitem[H]{H}
B Hopkins, {\it 30 Years of Bulgarian Solitaire}, Col. Math. J. {\bf 43} (2012), no. 2, 135--140.



%\bibitem[HS]{HS}
%Hopkins, B. and Sellers, J. A.
%{\em Exact Enumeration of Garden of Eden Partitions}, INTEGERS 7, no. 2 (2007), Article A19.

\end{thebibliography}
\end{document}